\documentclass[12pt]{article}
\renewcommand\b{{\bf b}}
\renewcommand\d{{\bf d}}

\newcommand\eps{\epsilon}
\newcommand\p{{\bf p}}

\renewcommand\r{{\bf r}}
\newcommand\s{{\bf s}}

\newcommand\x{{\bf x}}
\newcommand\y{{\bf y}}
\newcommand{\eref}[1]{$(\ref{#1})$}
\title{Some notes on applying computational divided differencing in 
optimization}
\author{Stephen A.~Vavasis\thanks{Department of Combinatorics and
Optimization, University of Waterloo, 200 University Ave.~W., Waterloo,
Ontario, Canada, N2L 3G1.  Supported in part by an NSERC Discovery
grant and a grant from the U.S.~Air Force Office of Scientific
Research.}}
\begin{document}
\maketitle
\begin{abstract}
We consider the problem of accurate computation of 
the finite difference $f(\x+\s)-f(\x)$
when $\Vert\s\Vert$ is very small.  
Direct evaluation of this difference
in floating point arithmetic
succumbs to cancellation error and yields 0 when $\s$ is sufficiently
small.  Nonetheless, accurate computation of this finite difference
is required by many optimization algorithms for a ``sufficient decrease''
test.  Reps and Rall proposed a programmatic transformation called
``computational divided differencing''  reminiscent of automatic
differentiation to compute these differences with high accuracy.  The
running time to compute the difference is a small constant multiple
of the running time to compute $f$.  Unlike automatic differentiation, however,
the technique is not fully general because of a difficulty with
branching code (i.e., `if' statements).  
We make several remarks about the application of computational
divided differencing to optimization.
One point is that the
technique can be used effectively as a stagnation test.
\end{abstract}

\section{Finite differences}

Many nonlinear optimization routines require a sufficient decrease
test on an iterate, which involves
computation of a finite difference.
For example the Armijo (also called ``backtrack'')
line-search requires testing an inequality of the form
$$f(\x+\alpha \p)-f(\x)\le \sigma \alpha \nabla f(\x)^T\p,$$
which is called a ``sufficient decrease'' condition.
The trust region method involves evaluating a ratio of the form
$$\frac{f(\x+\s)-f(\x)}{m(\x+\s)-m(\x)}$$
in which the numerator is a finite difference of the objective
while the denominator is a finite difference of a quadratic model.
See Nocedal and Wright \cite{NocedalWright} for more information about these
the Armijo line-search and trust region method.

A difficulty with these tests is that as the optimizer is approached,
the finite difference in $f$ is computed with decreasing accuracy
because the subtraction is dominated by cancellation.  This is because
 $\s$ becomes smaller so $f(\x+\s)$ and $f(\x)$ become closer.  This obstacle
is well known to implementors of optimization algorithms; see e.g.,
remarks by Nocedal and Wright on p.~62 or see 
Hager and Zhang \cite{HZ}.
The obvious workaround of using a Taylor approximation
in place of the finite difference (e.g., 
$f(\x+\s)-f(\x)\approx \nabla f(\x)^T\s$)
is not applicable in this setting
because the whole point of a sufficient decrease or ratio test is to
compare the finite difference to the Taylor approximation. 

Hager and Zhang \cite{HZ} observe that higher derivatives
can provide a sufficiently accurate approximation
in a sufficient-decrease test close to the root.
However, this solution is not completely
satisfactory because a high-order Taylor approximation rapidly loses
accuracy for larger $\s$, and it is not clear how to estimate the threshold
for switching between direct subtraction and Taylor approximation.

One might argue that in the case that the finite difference becomes so
inaccurate that the tests fail, the solution is nearly at hand anyway
so the algorithm can simply terminate.  This argument is valid in settings
when a low-accuracy solution is acceptable, but there are many settings
where high accuracy is desired.  One is the setting of testing
and comparing optimization routines.  In this case, one tries to drive
the algorithm into its asymptotic range to confirm predicted behavior.
A second setting is when the purpose of the optimization problem is not
to minimize the objective value per se but rather to drive the gradient
to zero.  This is the case in science and engineering problems in which
the objective is an energy potential, and the gradient of the potential
corresponds to a vector of unbalanced forces.  Driving the gradient to
zero is more important than minimizing the objective since a zero
gradient means a force equilibrium.

As observed and developed by Reps and Rall \cite{RepsRall},
the above-mentioned difficulty with finite differencing
is solved by rewriting the finite difference.
To take a simple example, suppose the objective function is the univariate
$f(x)=x^2$ at the point $x=1$.  If one evaluates $f(x+s)-f(x)$ for
$s=10^{-18}$ in IEEE double precision floating point arithmetic, 
then one obtains 0.  On the
other hand, if one rewrites $(x+s)^2-x^2$ as $2xs+s^2$, then the finite
difference can be computed to 16 significant digits of accuracy in IEEE
double-precision arithmetic for an $s$ of arbitrarily small magnitude.

For the rest of the paper, we use the notation $\mathcal{D}_f(\x,\s)$
to denote $f(\x+\s)-f(\x)$.
We now state the goal formally: we wish to compute $\mathcal{D}_f(\x,\s)$
with an error bounded above by $\Vert\s\Vert\eps_{\rm mach}c(f,\x)$, where
$\eps_{\rm mach}$ is machine-epsilon (unit roundoff) for the floating-point
arithmetic system, and $c(f,\x)$ specifies the level of roundoff error
present in the evaluation of $f(\x)$. 

This level of accuracy is attainable by the naive formula
$f(\x+\s)-f(\x)$ for large values of $\s$ since there is no
substantial cancellation between the terms in this case.  On the other
hand, it is clear that this level of accuracy cannot be attained for
small $\s$ using the computation $f(\x+\s)-f(\x)$; indeed, the error
is proportional to $\Vert \s\Vert$ in this case (since the answer will
be 0).  For extremely small values of $\s$, this level of accuracy is
attainable with the Taylor series approximation because
$f(\x+\s)-f(\x)-\nabla f(\x)^T\s$ is $O(\Vert\s\Vert^2)$ which will be
smaller than $\Vert\s\Vert\eps_{\rm mach}c(f,\x)$ for sufficiently
small $\s$.  However, for values of $\s$ that are small but not so
small (e.g., $\sqrt{\eps_{\rm mach}}$), neither of these straightforward
techniques is successful.

In this paper, we first briefly review computational divided
differencing as developed by Reps and Rall \cite{RepsRall},
who  show that the rewriting described
above can be carried out automatically on a program for computing $f$
by using a kind of chain rule to produce a new program that evaluates
the finite difference $\mathcal{D}_f(\x,\s)$ to high accuracy.  The
rules, which are reviewed in Section~\ref{sec:rules}, are similar to
automatic differentiation in the forward mode.  For a textbook
treatment of automatic differentiation as well as historical
background, see Chapter 8 of \cite{NocedalWright}.  The rewritten
program will have a running time at most a constant factor larger than
the original.  

Next, we review some issues with computational divided differencing
as applied to optimization.
The technique has major limitation not
present with automatic differentation, namely, it cannot handle the
general case of code branching, i.e., `if-else' statements.  An
explanation of the issue with branching and some partial workarounds
are presented in Section~\ref{sec:branching}.  

A second application of
computational divided differencing is in stagnation termination tests.
This application is discussed in Section~\ref{sec:stagnation}.

\section{Rules for computational divided differencing}
\label{sec:rules}

We assume we are given a program $P$ that takes as input $\x$ and produces
as output $f(\x)$.  The program consists of loops and assignment
statements.  Discussion of branching is postponed to
Section~\ref{sec:branching}.  To simplify the discussion, we assume
that each assignment statement in fact contains a single operation on
the right-hand side, either a floating-point binary or unary arithmetic
operation.  We further assume that the same variable does not occur
on both the left and right-hand sides of an assignment statement.
These assumptions are without loss of generality since 
complex assignment statements can be rewritten to obey these assumptions
via the introduction of intermediate variables.  

The output of the computational divided difference engine will be a new
program $P'$ that takes as input $\x$ and $\s$ produces as output both
$f(\x)$ and $\mathcal{D}_f(\x,\s)$.  It works by replacing every single
assignment statement of $P$, say one that assigns a variable $t$, with two
assignment statements in $P'$, one that computes $t$ and a second that
computes $\Delta t$, where $\Delta t$ is the difference between the value
$t$ would get when the program for
$f$ is invoked for $\x+\s$ minus the value computed for
$t$ when the program is invoked on $\x$.  This is done for all program
variables.  For the input variable $\x$, its finite difference
$\Delta\x$ is initialized to $\s$, that is, to the second input
of $P'$.
If a program variable $t$ is known to be a
{\em parameter}, that is, a number that does not depend on the
function argument $\x$, then one takes $\Delta t$ to be zero
since changing the input variable does not affect the parameter.

\begin{description}
\item[Addition \& subtraction.]
For `$+$' and `$-$', the finite differencing is linear.  In other
words, 
a statement of the form $t:=u+v$ in program $P$ causes
the following additional statement (i.e., in addition to
$t:=u+v$) to be inserted into $P'$:
$$\Delta t := \Delta u+\Delta v$$
and similarly for subtraction.

\item[Multiplication.]
For `$*$' (multiplication), the distributive law is used.
In particular,  if $P$ contains the statement
$t:=uv$, then $(u+\Delta u)(v+\Delta v)=
uv+v\Delta u + u\Delta v+\Delta u\cdot\Delta v$.
Thus, $P'$ contains the additional statement,
$$ \Delta t := u\cdot \Delta v + v\cdot \Delta u + 
\Delta u \cdot \Delta v.$$
The point is that the term $uv$ is ``precanceled'', meaning that
it is known to occur in both $uv$ and $(u+\Delta u)(v+\Delta v)$,
so it is simply omitted from the computation of $\Delta t$.

\item[Division.]
We break division
into two steps, reciprocation and multiplication.
Multiplication was covered above, and reciprocals are
covered below.

\item[Exponentiation.]
We rewrite $u^{v}$ as $\exp(v\log u)$ and use the rules
for log, exp and multiplication.  The rules for log and exp are below.
For the special case of squaring,
reciprocals, and square roots,
see below.  There are also special and more efficient
rules applicable for other
small constant rational-number exponents that are not presented
here.

\item[Reciprocals.]
For $t:=1/u$, the generated statement is
$$\Delta t:=\frac{-\Delta u}{u(u+\Delta u)}.$$

\item[Squaring.]
For $t:=u^2$, as mentioned in the introduction, we use
$$\Delta t:=2u\Delta u + (\Delta u)^2.$$

\item[Square roots.]
For $t:=\sqrt{u}$, we  use
$$\Delta t:= \frac{\Delta u}{\sqrt{u+\Delta u}+\sqrt{u}}.$$

\item[Exponential function.]
For $t:=\exp(u)$, we use
$$\Delta t := \exp(u)(\exp(\Delta u) - 1),$$
where the quantity in parentheses is evaluated as written
for $|\Delta u|\ge 1$ (since there is no risk of cancellation in
this case), or using a Taylor series (with the first term of 
`1' precanceled) if $|\Delta u|\le 1$.
We take 17 terms of the Taylor series, which is enough to
get full IEEE double precision.

\item[Logarithm.]
For $t:=\log(u)$, we use
$$\Delta t := \log(1+(\Delta u)/u).$$
The function \verb+log1p+, which is
available in several programming languages,
implements $\log(1+x)$ for $|x|$ small, and there are standard
polynomial approximations available for this function.
\end{description}

Many other standard math-library functions also have exact formulas
or good Taylor approximations for computational divided differencing.
A notable exception is the absolute value function, which is
considered in the next section.

\section{The difficulty with branching}
\label{sec:branching}

In this section we consider branching, that is, code with `if-else' blocks.
In the case of automatic differentiation, `if-else' blocks pose no difficulty
since, except for pathological situations, it is valid to move the
differentiation inside
the blocks.  But for computational divided differencing, `if-else' blocks
are problematic since it is possible that the evaluation of $f(\x+\s)$
would follow one branch while the evaluation of $f(\x)$ follows another,
so there is no obvious way to track the finite differences once
this divergence of execution paths occurs.

The first thing to observe is that if the code contains a branch
that causes a point of nondifferentiability, then it is not possible to 
compute accurate finite differences in the desired sense.  For example,
consider the computation of $f(x)=|x-\pi|$, where, as usual, $\pi=3.14\dots$.
The absolute value function implicitly involves an `if-else' block based
on the sign of its argument.  For the computation of $\mathcal{D}_f(x,s)$, if
$x$ is extremely close to $\pi$ and $s$ is very small, then
$\mathcal{D}_f(x,s)$ could be either $-s$ (if both $x$ and $x+s$
are less than $\pi$), $s$ (if both $x$ and $x+s$ are
greater than $\pi$) or some value in between.  Since there is no
way to distinguish these cases in floating point arithmetic if $x$ agrees
with $\pi$ to all digits of accuracy, 
we see that there is no way to guarantee
a successful evaluation of $\mathcal{D}_f(x,s)$.

Therefore, we limit attention to cases in which there is an `if-else'
block in which the resulting function is continuous and differentiable
at the breakpoint.  Even in this case, we do not know of a
general-purpose solution to handling `if-else' blocks.  However, we
have identified below four well-known examples of branching occurring
commonly in optimization that result in differentiable functions and a
method for handling these four cases.

\subsection{$\ell^2$ penalties}

A common differentiable function in optimization that involves
branching is the $\ell^2$ penalty function, that is, $p:=\max(0,x)^2$.
The finite difference for $p$ can be accurately evaluated as follows:
$$\Delta p:=\left\{
\begin{array}{ll}
2x\Delta x + \Delta x^2 & \mbox{if $x\ge 0$ and $x+\Delta x\ge 0$},\\
\max(0,x+\Delta x)^2-\max(0,x)^2 & \mbox{else.}
\end{array}\right.
$$
The point is that in the first case, the function is simply squaring
so we can use the squaring rule.  In the second case, one term or
the other is zero so there is no cancellation.

\subsection{Cubic splines}

Optimization objective functions sometimes use cubic splines to fit
data or to smooth a nonsmooth function.  Suppose the $C^2$ cubic spline 
$s(x)$ is
presented as follows.  There is a series of knots
$\xi_1<\xi_2<\cdots<\xi_k$.  Within each subinterval
$[\xi_i,\xi_{i+1}]$, there are two representations of the spline
function, one of the form
$a_i(x-\xi_i)^3+b_i(x-\xi_i)^2+c_i(x-\xi_i)+d_i$ and the other of the
form $m_i(x-\xi_{i+1})^3+n_i(x-\xi_{i+1})^2+o_i(x-\xi_{i+1})+p_i$, and
these should agree to all significant decimal places.  Furthermore,
it should be the case that $p_i=d_{i+1}$ (to all decimal digits) for
each $i$ to ensure continuity across breakpoints.  There are additional
conditions to ensure $C^2$ continuity that we do not
specify here.  Finally, we assume that for
$x<\xi_1$ there is a specification of a cubic $m_0(x-\xi_1)^3+\cdots+p_0$,
and for $x>\xi_k$ there is a specification $a_k(x-\xi_k)^3+\cdots+d_k$.

Then the finite difference routine to evaluate $\mathcal{D}_s(x,\Delta x)$
for the
call $t:=s(x)$ is as follows.  Without loss of generality, $\Delta x\ge 0$
since otherwise we can exchange the roles of $x$ and $x+\Delta x$ and
then invert the sign of $\Delta t$ at the end of the computation.

If $x$ and $x+\Delta x$ are in the
same subinterval, then we use a finite difference formula for
cubics (easily obtained by binomially expanding each power in
$(x+\Delta x)$ and then subtracting like terms).

If $x$ and $x+\Delta x$ are in different intervals, say
$x\in[\xi_i,\xi_{i+1}]$ and $x+\Delta x\in[\xi_j,\xi_{j+1}]$ where
$j>i$ and where we identify $\xi_0=-\infty$ and $\xi_{k+1}=\infty$, then
we evaluate $\Delta t$ using the telescoping formula:
$$\Delta t= (s(x+\Delta x)- s(\xi_j))+\sum_{l=i+1}^{j-1} (s(\xi_{l+1})-s(\xi_l))
+(s(\xi_i)-s(x)).$$
The first and last terms are evaluated via finite differencing
of a polynomial, noting
that the constant terms precancel
(because the formula for $s$ between knots includes
the value at the knot as the constant coefficient).  The middle terms
are evaluated by subtracting the relevant constant coefficients.

\subsection{Termination tests for iterative methods}

If the objective function includes an iterative loop that terminates
when a tolerance is sufficiently small, there is apparently no easy
way to handle this case using the framework explained herein.  A workaround
is to determine the maximum number of iterations needed
by the iteration taken over all data within the
feasible region, and then replacing the conditional loop termination with
a simple loop that always runs up to the maximum number of iterations.

\subsection{Gaussian elimination}

Gaussian elimination with partial or complete pivoting 
to solve a square
system of linear equations involves branching
on which row or column is selected for the pivot.  Nonetheless, the
overall computation of $\x:=A^{-1}\b$ is amenable to finite differencing,
and finite differences can be evaluated via:
\begin{eqnarray*}
\Delta\x&=&(A+\Delta A)^{-1}(\b+\Delta\b)-A^{-1}\b \\
&=& A^{-1}([(I+\Delta A\cdot A^{-1})^{-1}-I]\b+(I+\Delta A\cdot A^{-1})^{-1}
\Delta \b).
\end{eqnarray*}
The factor in square brackets can be evaluated by a Taylor expansion,
precanceling $I$, if 
$\Vert \Delta A\cdot A^{-1}\Vert\le 1/2$ in an induced matrix norm
(e.g., the matrix $\infty$-norm), else it can be evaluated by direct
subtraction.  The second term can be evaluated by direct
expansion.

\section{Stagnation termination test}
\label{sec:stagnation}

Assume for this section that an descent-based optimization algorithm
is under consideration, that is, one in which an objective function
(or perhaps a penalized objective function in the presence of
constraints) decreases from one iterate to the next.  Most such
optimization algorithms involve a stagnation test that is triggered
when insufficient progress is made for some number of consecutive
iterations.

The question is how to measure insufficient progress.  The obvious way
is to check the decrease in the objective function.  For example, a
stagnation test could be triggered if there is no relative decrease of
more, than, e.g., $10^{-15}$ in the objective function for three
successive iterations.

The problem with this test is that it could be triggered too early in
the case that the objective function is a sum of terms, some of which
are so large that the smaller terms are lost in the least-significant digits
of the objective value.  This can happen in a potential-energy formulation
mentioned in the introduction when part of the configuration is
minimized, and the energy contribution from the part already minimized
is a large number in absolute terms, but another part is still
rapidly evolving.

Computational divided differencing provides an alternative and perhaps more
robust method to determine stagnation.  The stagnation test is based
on the following observation.  In exact arithmetic, if $\x_1,\x_2,\x_3$
are successive iterates, then obviously
$$f(\x_1)-f(\x_3)=(f(\x_1)-f(\x_2))+(f(\x_2)-f(\x_3)),$$
i.e.,
$$\mathcal{D}_f(\x_3,\x_1-\x_3)
=
\mathcal{D}_f(\x_2,\x_1-\x_2)+
\mathcal{D}_f(\x_3,\x_2-\x_3).$$
If the three finite differences are evaluated using the above rules,
we would again expect approximate equality
to hold.  Our proposed stagnation test is when the left-hand side
is much smaller than the right-hand side, say by a factor of 2.
This means that the progress predicted by exact finite differences
is not observed, so progress is no longer possible.  The stagnation test
can also be applied over a longer sequence of steps.

To explain in more detail why this method may be more suitable than
directly using the objective function, 
consider an objective function of the form
$f(\x)=\frac{1}{2}\x^TM\x+\d^T\x$ where $M$ is symmetric positive definite.
The minimizer is clearly at $\x^*=-M^{-1}\d$, and the minimum
objective value is $-\d^TM^{-1}\d/2$.
Suppose that it is possible to evaluate both $\x^TM\x$ and
$M\x$ in a
forward-accurate sense for any $\x$, i.e., the computed $\hat\y\approx M\x$ 
and $\hat\alpha\approx\x^TM\x$ in the
presence of roundoff error satisfy 
$\Vert\y-\hat\y\Vert\le
c\eps_{\rm mach}\Vert\y\Vert$  and 
$|\alpha -\hat\alpha |\le
c\eps_{\rm mach}\alpha$
for a small $c$.  
This is possible if $M$ is well-conditioned.
For a general unstructured
ill-conditioned matrix, this bound is not possible (consider the case
that $\x$ is close to the eigenvector of the smallest eigenvalue, in which
case the computation of $M\x$ is likely to be inaccurate), but
in the case that $M$ has an appropriate
partially separable representation (e.g.,
$M$ is diagonal), such a bound holds even when $M$ is ill-conditioned.

Now consider an iterate $\x$ that is of the form $\x_1=-M^{-1}\d+\r$,
where $\r$ is the error vector and assumed to be generic.  Note that 
\begin{eqnarray*}
f(\x_1)&=&\frac{1}{2}\x_1^TM\x_1+\d^T\x_1 \\
&=& \frac{1}{2}(M^{-1}\d-\r)^TM(M^{-1}\d-\r)-\d^T(M^{-1}\d-\r) \\
&=& -\frac{1}{2}\d^TM^{-1}\d+\r^TM\r
\end{eqnarray*}
A stagnation test based on the objective value
will determine that no progress is possible
once the second term on the right-hand side is much smaller
than the first.
Thus, if 
\begin{equation}
\r^TM\r \le \d^TM^{-1}\d \eps_{\rm mach},
\label{eq:termcrit1}
\end{equation}
no further progress is possible 
for reducing the
objective function in this floating point arithmetic.

Now consider the proposed stagnation test.  Note that the finite
difference function
\begin{equation}
\mathcal{D}_f(\x,\Delta\x)=
(M\x+\d)^T\Delta\x + \frac{1}{2}\Delta\x^TM\Delta\x,
\label{eq:mathcaldf}
\end{equation}
involves products with $M$ that we are assuming are computed accurately.
Suppose $\x_i=-M^{-1}\d+\r_i$ for $i=1,\ldots,3$.  Then for $i=1,2$,
\begin{equation}
\mathcal{D}_f(\x_i,\x_{i+1}-\x_i)=(\r_{i+1}-\r_i)M\r_i + (\r_{i+1}-\r_i)^TM
(\r_{i+1}-\r_i)/2
\label{eq:mathcaldf2}
\end{equation}
Under the assumption discussed earlier plus some further assumptions
(for example, $|\r_{i+1}^TM\r_i| \le c \r_i^TM\r_i$, for a $c<1$ as would be
true for a quasi-Newton method), these differences are highly accurately
computed assuming that the $\r_i$'s are known accurately.

In fact, the inaccuracy in the test will arise primarily because
the $\r_i$'s are not known accurately.
If $\x_i=-M^{-1}\d+\r_i$, then
$M\x_i+\d=M(-M^{-1}\d+\r_i)+\d=-M\cdot M^{-1}\d+\d+M\r_i=M\r_i$.  The
sum $-\d+\d$ cannot be precanceled (because the algorithm does not
represent iterates explicitly in the form $-M^{-1}\d+\r_i$), so
there will be cancellation error in evaluating $M\x_i+\d$,
In other words, the program will
make substantial errors in the evaluation of \eref{eq:mathcaldf2}
once $\r_i$ is sufficiently small with respect to $M^{-1}\d$.  This occurs
when 
\begin{equation}
\Vert \r_i\Vert \approx c\Vert M^{-1}\d\Vert\cdot\eps_{\rm mach}
\label{eq:termcrit2}
\end{equation}
for a small constant $c$.
Notice the difference between \eref{eq:termcrit1}
and \eref{eq:termcrit2}.  A sufficient condition in terms of norms
to imply \eref{eq:termcrit1} is
$$\frac{\Vert \r\Vert\cdot\Vert M\Vert }{\Vert\d\Vert}\le c\sqrt{\eps_{\rm mach}}$$
whereas a sufficient condition to imply \eref{eq:termcrit2} is
$$\frac{\Vert \r\Vert\cdot\Vert M\Vert }{\Vert\d\Vert}\le c\eps_{\rm mach}.$$
Thus, the algorithm is able to make much
more progress (and indeed, obtain the solution $-M^{-1}\d$ to full machine
precision) with the second stagnation test rather then the first.

\section{Acknowledgment}
The author acknowledges help from Chris Bischof, who led him to the paper
by Reps and Rall.

\section{Conclusions}
\label{sec:conclusions}

We have made several remarks about automatically 
computing finite differences using Reps and Ralls' computational
divided differencing scheme.  This scheme, when
given a program to compute $f(\x)$, generates
a new program that computes the finite difference
$f(\x+\s)-f(\x)$ to high relative accuracy even when $\s$ is small.
These finite differences are useful for sufficient-decrease
tests and ratio tests in optimization algorithms.
The technique is not fully general due to the problem with `if-else'
statements but is applicable to many optimization problems in which
the objective function is a sequence of complicated arithmetic expressions
such as energy functionals in computational mechanics.
We have also suggested a test based on computational
divided differencing may be used to detect stagnation
in optimization routines more accurately than merely using the
objective function.
\bibliography{icov}

\begin{thebibliography}{1}

\bibitem{HZ}
W.~Hager and H.~Zhang.
\newblock A new conjugate gradient method with guaranteed descent and an
  efficient line search.
\newblock {\em SIAM J. Optimiz.}, 16(1):170--192, 2005.

\bibitem{NocedalWright}
J.~Nocedal and S.~Wright.
\newblock {\em Numerical Optimization, 2nd Edition}.
\newblock Springer, New York, 2006.

\bibitem{RepsRall}
T.~W. Reps and L.~B. Rall.
\newblock Computational divided differencing and divided-difference
  arithmetics.
\newblock {\em Higher-order and symbolic computation}, 16:93--149, 2003.

\end{thebibliography}
\bibliographystyle{plain}
\end{document}